\documentstyle[12pt]{article}
\begin{document}
\newcommand{\text}[1]{\mbox{{\rm #1}}}  
\newcommand{\gd}{\delta}
\newcommand{\itms}[1]{\item[[#1]]} 
\newcommand{\nin}{\in\!\!\!\!\!/}
\newcommand{\sub}{\subset} 
\newcommand{\cntd}{\subseteq}     
\newcommand{\go}{\omega} 
\newcommand{\Pa}{P_{a^\nu,1}(U)} 
\newcommand{\fx}{f(x)}  
\newcommand{\fy}{f(y)} 
\newcommand{\gD}{\Delta}
\newcommand{\gl}{\lambda} 
\newcommand{\gL}{\Lambda} 
\newcommand{\half}{\frac{1}{2}} 
\newcommand{\sto}[1]{#1^{(1)}}
\newcommand{\stt}[1]{#1^{(2)}}
\newcommand{\Z}{\hbox{\sf Z\kern-0.720em\hbox{ Z}}}
\newcommand{\singcolb}[2]{\left(\begin{array}{c}#1\\#2
\end{array}\right)} 
\newcommand{\ga}{\alpha}
\newcommand{\gb}{\beta} 
\newcommand{\gga}{\gamma}
\newcommand{\ul}{\underline} 
\newcommand{\ol}{\overline} 
\newcommand{\qed}{\kern 5pt\vrule height8pt width6.5pt depth2pt}
\newcommand{\Lrraro}{\Longrightarrow}
\newcommand{\Nb}{|\!\!/}
\newcommand{\NN}{{\rm I\!N}}
\newcommand{\bsl}{\backslash}     
\newcommand{\gt}{\theta}
\newcommand{\op}{\oplus}
\newcommand{\Op}{\bigoplus}          
\newcommand{\CR}{{\cal R}}
\newcommand{\tr}{\bigtriangleup}
\newcommand{\grr}{\omega_1} 
\newcommand{\ben}{\begin{enumerate}}
\newcommand{\een}{\end{enumerate}}
\newcommand{\ndiv}{\not\mid}
\newcommand{\bab}{\bowtie}
\newcommand{\hal}{\leftharpoonup}
\newcommand{\har}{\rightharpoonup}
\newcommand{\ot}{\otimes}
\newcommand{\OT}{\bigotimes}
\newcommand{\bwe}{\bigwedge}
\newcommand{\gep}{\varepsilon}
\newcommand{\gs}{\sigma} 
\newcommand{\rbraces}[1]{\left( #1 \right)}
\newcommand{\bbox}{$\;\;\rule{2mm}{2mm}$}
\newcommand{\sbraces}[1]{\left[ #1 \right]}
\newcommand{\bbraces}[1]{\left\{ #1 \right\}}
\newcommand{\OO}{_{(1)}}
\newcommand{\TT}{_{(2)}}
\newcommand{\FF}{_{(3)}}
\newcommand{\minus}{^{-1}}
\newcommand{\CV}{\cal V} 
\newcommand{\CVs}{\cal{V}_s} 
\newcommand{\un}{U_q(sl_n)'}
\newcommand{\on}{O_q(SL_n)'}
\newcommand{\slq}{U_q(sl_2)}
\newcommand{\olq}{O_q(SL_2)}
\newcommand{\UU}{U_{(N,\nu,\go)}}
\newcommand{\HH}{H_{n,q,N,\nu}} 
\newcommand{\ct}{\centerline}
\newcommand{\bs}{\bigskip}
\newcommand{\qua}{\rm quasitriangular}   
\newcommand{\ms}{\medskip}
\newcommand{\noin}{\noindent}
\newcommand{\mat}[1]{$\;{#1}\;$}
\newcommand{\raro}{\rightarrow}
\newcommand{\map}[3]{{#1}\::\:{#2}\raro{#3}}
\newcommand{\alg}{{\rm Alg}}
\def\newtheorems{\newtheorem{theorem}{Theorem}
                 \newtheorem{lemma}[theorem]{Lemma}
                 \newtheorem{Theorem}{Theorem}
                 \newtheorem{Lemma}[Theorem]{Lemma}}
\newtheorems
\newcommand{\proof}{\par\noindent{\bf Proof:}\quad}
\newcommand{\dmatr}[2]{\left(\begin{array}{c}{#1}\\
                            {#2}\end{array}\right)}
\newcommand{\doubcolb}[4]{\left(\begin{array}{cc}#1&#2\\
#3&#4\end{array}\right)}
\newcommand{\qmatrl}[4]{\left(\begin{array}{ll}{#1}&{#2}\\
                            {#3}&{#4}\end{array}\right)}
\newcommand{\qmatrc}[4]{\left(\begin{array}{cc}{#1}&{#2}\\
                            {#3}&{#4}\end{array}\right)}
\newcommand{\qmatrr}[4]{\left(\begin{array}{rr}{#1}&{#2}\\
                            {#3}&{#4}\end{array}\right)}
\newcommand{\smatr}[2]{\left(\begin{array}{c}{#1}\\
                            \vdots\\{#2}\end{array}\right)}

\newcommand{\ddet}[2]{\left[\begin{array}{c}{#1}\\
                           {#2}\end{array}\right]}
\newcommand{\qdetl}[4]{\left[\begin{array}{ll}{#1}&{#2}\\
                           {#3}&{#4}\end{array}\right]}
\newcommand{\qdetc}[4]{\left[\begin{array}{cc}{#1}&{#2}\\
                           {#3}&{#4}\end{array}\right]}
\newcommand{\qdetr}[4]{\left[\begin{array}{rr}{#1}&{#2}\\
                           {#3}&{#4}\end{array}\right]}

\newcommand{\qbracl}[4]{\left\{\begin{array}{ll}{#1}&{#2}\\
                           {#3}&{#4}\end{array}\right.}
\newcommand{\qbracr}[4]{\left.\begin{array}{ll}{#1}&{#2}\\
                           {#3}&{#4}\end{array}\right\}}

\title{Semisimple Hopf Algebras of Dimension $pq$ are Trivial}
\author{Pavel Etingof and Shlomo Gelaki 
\\Department of Mathematics\\
Harvard University\\Cambridge, MA 02138}
\date{December 9, 1997}
\maketitle
This paper makes a contribution to the problem of classifying 
finite-dimensional semisimple Hopf algebras $H$ over an algebraically 
closed field $k$ of characteristic $0.$ Specifically, we show that if $H$ 
has dimension $pq$ for primes $p$ and $q$ then $H$ is trivial, that is, $H$ 
is either a group algebra or the dual of a group algebra.

Previously known cases include dimension 
$2p$ \cite{m1}, dimension $p^2$ \cite{m2}, and dimensions $3p,$ $5p$ 
and $7p$ \cite{gw}. Westreich and the second author also obtained the 
same result for $H$ which is, along with 
its dual $H^*,$ of Frobenius type (i.e. the dimensions of their 
irreducible representations divide the dimension of $H$) [GW, Theorem 3.5]. 
They concluded with the conjecture that any semisimple Hopf algebra $H$ of 
dimension $pq$ over $k$ is trivial. 

In this paper we use Theorem 1.4 in \cite{eg} to prove that both $H$ and
$H^*$ are of Frobenius type, and hence prove this conjecture.

Throughout $k$ will denote
an algebraically closed field of characteristic $0,$ and 
$p$ and $q$ will denote two prime numbers satisfying $p<q.$ 
We let $G(H)$ denote the group of grouplike 
elements of a Hopf algebra $H,$ and $D(H)$ denote its Drinfel'd double. 
Recall that $H$ and $H^{*cop}$ are Hopf subalgebras of $D(H),$ and let
$i_{H}:H\hookrightarrow D(H)$ and $i_{H^{*cop}}:H^{*cop}\hookrightarrow D(H)$
denote the corresponding inclusion maps. 

The following result will prove very useful in the sequel.
\begin{Theorem} {\bf [GW, Theorem 2.1]}\label{sara}
Let $H$ be a semisimple Hopf algebra
of dimension $pq,$ where $p<q$ are two prime numbers. Then:
\ben
\item $|G(H)|\ne q.$
\item If $|G(H)|=p$ then $q=1(mod\,p).$
\een
\end{Theorem}

Suppose there exists a non-trivial semisimple Hopf algebra $H$ of 
dimension $pq.$ Then by Theorem \ref{sara}, $G(H)$ and $G(H^*)$ have either 
order $1$ or $p.$ By \cite{r}, 
$G(D(H))=G(H^*)\times G(H)$ has therefore either order $1,p$ or $p^2.$ Since 
by \cite{r}, the order of $G(D(H)^*)$ divides the order of $G(D(H)),$ it 
follows that $G(D(H)^*)$ has either order $1,p$ or $p^2$ too. In the 
following we show that this leads to a contradiction, and hence 
conclude that $H$ is trivial.
\begin{Lemma}\label{start}
Let $H$ be a semisimple Hopf algebra of dimension $pq$ over $k,$ and
let $V,U$ be $p-$dimensional irreducible representations of $D(H).$ If
$V\ot U$ contains a $1-$dimensional representation $\chi,$ then it must 
contain another $1-$dimensional representation $\chi'\ne\chi.$ 
\end{Lemma}
\proof Without loss of generality we can assume that $\chi$ is trivial and 
hence that $U=V^*.$ Otherwise we can replace $U$ with $U\ot \chi^{-1}.$

Suppose on the contrary that $V\ot V^*$ does not contain a non-trivial
$1-$dimensional representation. By \cite{eg},
the dimension of any irreducible representation of $D(H)$ divides $pq,$ hence
$V\ot V^*$ is a 
direct sum of the trivial representation $k,$ $p-$dimensional irreducible 
representations of $D(H)$ 
and $q-$dimensional irreducible representations of $D(H).$ Therefore we
have that $p^2=1+ap+bq.$ Clearly $b>0.$ Let $W$ be a $q-$dimensional 
irreducible 
representation of $D(H)$ such that $W\subset V\ot V^*.$ Since $0\ne 
Hom_{D(H)}(V\ot 
V^*,W)=Hom_{D(H)}(V,W\ot V)$ we have that $V\subset W\ot V.$ Since 
$W\ot V$ has no linear constituent (because $dimV\ne dim W$), $dim(W\ot V)=pq$ 
and $W\ot V$ contains a 
$p-$dimensional irreducible representation of $D(H)$ it follows that 
$W\ot V=V_1\oplus\cdots \oplus V_q$ where $V_i$ is a $p-$dimensional 
irreducible representation of $D(H)$ with $V_1=V.$ 

We wish to show that 
for any $i=1,\dots,q$ there exists a $1-$dimensional representation 
$\chi_i$ such that $V_i=V\ot \chi_i.$  
Suppose on the contrary that this is not true for some $i.$ Then $V\ot 
V_i^*$  has no linear constituent, 
hence it must be a direct sum of $p$ $p-$dimensional irreducible 
representations of 
$D(H).$ Therefore, $W\ot (V\ot V_i^*)$ has no linear constituent.
But, $(W\ot V)\ot V_i^*=(V_i\oplus\cdots )\ot V_i^*=k\oplus\cdots $ 
which is a contradiction. 
Therefore, $W\ot V=V\ot (\chi_1\oplus\cdots\oplus\chi_q).$ 
Note that for all $i,$ $\chi_i$ is uniquely determined. Indeed, if there 
exists $\chi_i'\ne \chi_i$ such that $V\ot \chi_i=V\ot \chi_i',$ 
then $\chi_i'\chi_i^{-1}$ is a non-trivial linear constituent of $V\ot V^*$ 
which is a contradiction. 

We now consider two possible cases. First, suppose that any $1-$dimensional 
representation $\chi$ such that $\chi\ot W=W,$ is trivial. Then for any 
$q-$dimensional irreducible representation $W'$ of $D(H),$ $W\ot W'$ 
contains no more than one $1-$dimensional representation. Therefore, since 
$V\ot V^*$ contains $b\ge 1$ $q-$dimensional irreducible representations of 
$D(H),$ it follows that $W\ot V\ot V^*$ contains at most $b$ $1-$dimensional 
representations. Therefore, $b\ge q,$ and $p^2=dim(V\ot V^*)\ge bq\ge q^2,$ 
which is a contradiction. 

Second, suppose there exists a non-trivial $1-$dimensional representation 
$\chi$ such that $\chi\ot W=W.$ Then $\chi\ot W\ot V=W\ot V$ which is 
equivalent to $\chi\ot V\ot (\chi_1\oplus\cdots\oplus\chi_q)=V\ot 
(\chi_1\oplus\cdots\oplus\chi_q).$ Write 
$\chi_1\oplus\cdots\oplus\chi_q=\sum _{\ga \in G(D(H)^*)}m_{\ga}\ga.$ 
Then $V\ot \sum _{\ga \in G(D(H)^*)}m_{\ga \chi^{-1}}\ga=
V\ot \sum _{\ga \in G(D(H)^*)}m_{\ga}\ga \chi=
V\ot \sum _{\ga \in G(D(H)^*)}m_{\ga}\ga.$ But, $V\ot \ga_1=V\ot \ga_2$ 
implies $\ga_1=\ga_2,$ so $m_{\ga}=m_{\ga \chi^{-1}}.$ Since the order of 
$\chi$ in $G(D(H)^*)$ is divisible by $p,$ this implies 
that $q=\sum_{\ga}m_{\ga}$ is divisible by $p$ which is a contradiction.
This concludes the proof of the lemma.\qed 
\begin{Lemma}\label{nt}
Let $H$ be a semisimple Hopf algebra of dimension $pq$ over $k.$ Then  
$G(D(H)^*)$ is non-trivial.
\end{Lemma} 
\proof 
Suppose on the contrary that $|G(D(H)^*)|=1.$ By \cite{eg}, the 
dimension of any irreducible representation of $D(H)$ divides $pq,$ hence 
$p^2q^2=1+ap^2+bq^2$ for some integers $a,b>0.$ In particular $D(H)$ has 
a $p-$dimensional irreducible representation $V.$ But then by Lemma 
\ref{start}, $D(H)$ must have a non-trivial $1-$dimensional representation 
which is a contradiction.\qed
\begin{Lemma}\label{npsq}
Let $H$ be a non-trivial semisimple Hopf algebra of dimension $pq$
over $k.$ Then 
the order of $G(D(H)^*)$ is not $p^2.$
\end{Lemma}
\proof Suppose on the contrary that $|G(D(H)^*)|=p^2.$ By \cite{r}, any 
$\chi\in 
G(D(H)^*)$ is of the form $\chi=g\ot \ga ,$ where $g\in G(H)$ and $\ga 
\in G(H^*).$ In particular, either $G(H)$ or $G(H^*)$ is non-trivial, and
hence by Theorem \ref{sara}, either $G(H)$ or $G(H^*)$ has order $p.$ 
Therefore we must have that both $G(H)$ and $G(H^*)$ are of order $p,$
and hence that $G(D(H)^*)=G(H)\times G(H^*).$ This implies that there 
exists a 
non-trivial $g\in G(H)$ so that $g\ot 1\in G(D(H)^*).$ By \cite{r}, $1\ot 
g\in G(D(H))$ is central in $D(H).$ In particular $g$ is central in $H,$ 
and hence $H$ 
possesses a normal Hopf subalgebra isomorphic to $kC_p,$ which implies 
that $H$ is commutative (see the proof of Theorem 2.1 in \cite{gw}). This
is a contradiction and the result follows.\qed
\begin{Lemma}\label{np}
Let $H$ be a non-trivial semisimple Hopf algebra of dimension $pq$ over
$k.$  Then there exists a quotient Hopf algebra 
$A$ of $D(H)$ of dimension $pq^2$ which contains $H$ and $H^{*cop}$ as 
Hopf subalgebras, and $H\cap H^{*cop}=kG(A)=kC_p.$ 
\end{Lemma}
\proof
By Lemmas \ref{nt} and \ref{npsq}, the order of $G(D(H)^*)$ is $p.$
If either $G(H)$ or $G(H^*)$ is trivial then either $H^*$ 
or $H$ respectively contains a central grouplike element. Hence, 
as in the proof of Lemma \ref{npsq}, $H$ is trivial. Therefore, $G(H)$ and 
$G(H^*)$ have both order $p.$ 

Since $D(H)$ is of Frobenius type \cite{eg} and since by Theorem 
\ref{sara}, $q=1(mod\,p),$ we must have that 
$p^2q^2=p+ap^2+bq^2,$ where $a=\frac{q^2-1}{p}$ and $b=p(p-1).$ 
In particular $D(H)$ has irreducible representations of dimension $p$ and 
$q.$ We wish to show that if $V$ and $U$ are two irreducible representations 
of $D(H)$ of dimension $p$ then $V\ot U$ is a direct sum of $1-$dimensional 
irreducible representations of $D(H)$ and $p-$dimensional irreducible 
representations of $D(H)$ only. Indeed, 
by Lemma \ref{start}, either $V\ot U$ does not contain any $1-$dimensional
representation or it must contain at least two different $1-$dimensional 
representations. But if it contains two different $1-$dimensional 
representations, then since the $1-$dimensional representations of $D(H)$  
form a cyclic group of order $p,$ it follows that $V\ot U$ contains all
the $p$ $1-$dimensional 
representations of $D(H).$ We conclude that either $p^2=ap+bq$ or 
$p^2=p+ap+bq.$ At any rate $b=0,$ and the result follows. 
Therefore, the 
subcategory ${\cal C}$ of $Rep(H)$ generated by the $1$ and $p-$dimensional 
representations of $D(H)$ gives rise to a quotient Hopf algebra $A$ of 
$D(H)$ of 
dimension $pq^2$ whose irreducible representations have either dimension 
$1$ or $p$ ( $A$ is the quotient of $D(H)$ over the Hopf ideal 
$I=\sum_{V\not\in{\cal C}}End_k(V)$). Let $\phi$ denote the corresponding 
surjective homomorphism, and consider the sequences:
$$H\stackrel{i_H}{\hookrightarrow}D(H)\stackrel{\phi}{\raro}A\;\;\;and\;\;\;
H^{*cop}\stackrel{i_{H^{*cop}}}{\hookrightarrow}D(H)\stackrel{\phi}{\raro}A.$$ 
Since $\phi(H),$ $\phi(H^{*cop})$ are Hopf subalgebras of $A$ it follows 
that they have dimension $1,$ $p,$ $q$ or $pq$ \cite{nz}. If, say, 
$dim(\phi(H))=q$ then $\phi(H)=kC_q$ is a group algebra \cite{z}. Since 
$(kC_q)^*\cong kC_q$ the surjection of Hopf algebras $\phi:H\raro kC_q$ 
gives rise to an inclusion of Hopf algebras $\phi^*:kC_q\raro H^*.$ 
Since $\phi^*(kC_q)\subset kG(H^*)$ it follows that  
$q$ divides $|G(H^*)|=p$ which is impossible. Similarly, 
$dim(\phi(H^{*cop}))\ne q,$ and we 
conclude that $\phi(H),$ $\phi(H^{*cop})$ have dimension $1,$ $p$ or $pq.$ 
Since $A=\phi(H)\phi(H^{*cop})$ is of dimension $pq^2$ it follows that the 
dimensions must equal $pq$ 
and hence that $H$ and $H^{*cop}$ can be considered as Hopf subalgebras of 
$A.$ Furthermore, $G(H)$ and $G(H^{*cop})$ are mapped onto $G(A),$ 
otherwise $|G(A)|$ is divisible by $p^2$ which is impossible. Therefore 
$G(H)\cap G(H^{*cop})=G(A)$ in 
$A.$ Since $H\ne H^{*cop}$ in $A,$ this implies that $H\cap 
H^{*cop}=kG(A)=kC_p.$\qed

We are ready to prove our main result.
\begin{Theorem}\label{main}
Let $H$ be a semisimple Hopf algebra of dimension $pq$ over $k,$ where $p$ 
and $q$ are distinct prime numbers. Then $H$ is trivial.
\end{Theorem}
\proof Suppose on the contrary that $H$ is non-trivial. We wish to prove that
$H$ and $H^*$ are of Frobenius type. Indeed,
consider the Hopf algebra quotient $A$ which exists by Lemma \ref{np}.
For any $1$-dimensional representation $\chi$ of $A$ we define the 
induced representations 
$$V_+^\chi=A\ot _H \chi\;\;\;and\;\;\;V_-^\chi=A\ot _{H^*}\chi.$$
Note that since finite-dimensional Hopf algebras are free
over any of their Hopf subalgebras \cite{nz}, it follows that
$V_+^\chi=H^*\ot _{H\cap H^*}\chi=H^*\ot
_{kC_p}\chi,$ and $V_-^\chi=H\ot _{H\cap H^*}\chi=H\ot _{kC_p}\chi.$

We first show that if $\chi$ is non-trivial then $\chi$ is non-trivial on 
$H$ and $H^*$ as well. Indeed, for any non-trivial 
$1-$dimensional representation $\tilde\chi$ we have by Frobenius 
reciprocity that $Hom_A(V_+^\chi,\tilde\chi)=Hom_H(\chi,\tilde\chi).$ 
Therefore, if $\chi$ is trivial on $H$ then $\tilde\chi$ is also trivial 
on $H$ (since $\tilde\chi=\chi^l$ as $\chi$ generates the group of
$1-$dimensional representations of $A$), and 
$Hom_A(V_+^\chi,\tilde\chi)=k.$ But then $V_+^\chi$ is a sum of $p$ 
$1-$dimensional representations and $p-$dimensional irreducible 
representations which contradicts the fact that $p$ does not divide $q.$
Similarly, $\chi$ is non-trivial on $H^*.$

Next, since $dim(Hom_A(V_+^\chi,\tilde\chi))=dim(Hom_H(\chi,\tilde\chi))=
\gd_{\chi,\tilde\chi}$ for any $1-$dimensional representation
$\tilde\chi$ and $dimV_+^\chi=q,$ it follows that
$$V_+^\chi=\chi\oplus V_1\oplus\cdots\oplus 
V_{\frac{q-1}{p}}$$ 
where $V_i$ is a $p-$dimensional irreducible 
representation of $A$ for all $1\le i\le \frac{q-1}{p}.$ 
We wish to show that the $V_i$'s do not depend on $\chi.$ Indeed, let
$\chi$ be a non-trivial $1-$dimensional representation of $A.$ Then
any $1-$dimensional representation of $A$ is of the form $\chi^l$ for 
some $0\le l\le p-1,$ and 
$V_+^{\chi^l}=V_+^\chi\ot \chi^{l-1}=\chi^l\oplus
(V_1\ot \chi^{l-1})\oplus\cdots\oplus (V_{\frac{q-1}{p}}\ot \chi^{l-1}).$
But by Lemma \ref{start}, $\chi^{l-1}\subset V_i\ot V_i^*$ for all $i,$
and hence $V_i\ot \chi^{l-1}\cong V_i$ which proves the claim.
Similarly, 
$V_-^\chi=\chi\oplus W_1\oplus\cdots\oplus W_{\frac{q-1}{p}}$
where $W_i$ is a $p-$dimensional irreducible representation of $A$
for all $1\le i\le \frac{q-1}{p}.$ 

Finally, we wish to show that the set
$$C=\{\chi_{|H^*}|\chi\;{\rm is}\;{\rm a}\;1-{\rm
dimensional}\;{\rm representation}\;{\rm of}\;A\}\cup 
\{V_{i|H^*}|1\le i\le (q-1)/p\}$$ 
is a full set of 
representatives of 
the irreducible representations of $H^*.$ Indeed, consider the regular 
representation of $H^*:$ 
\begin{eqnarray*}
H^*=H^*\ot _{kC_p}kC_p=\bigoplus_{\chi}H^*\ot 
_{kC_p}\chi=\bigoplus_{\chi}V_{+|H^*}^{\chi}=\bigoplus_{\chi}\chi
\oplus\bigoplus_i pV_{i|H^*}. 
\end{eqnarray*}
Since in the regular representation any
representation can occur no more times than its dimension we get that 
$V_{i|H^*}$
are irreducible. By summation of degrees, $C$ is a full set
of representatives, and since $dim(V_{i|H^*})=p$ we have that $H^*$ is of
Frobenius type. Similarly $H$ is of Frobenius type. 
But Theorem 3.5 in \cite{gw} 
states that if this is the case then $H$ is trivial. This concludes the 
proof of the theorem.\qed

\noin
{\bf Acknowledgments}
The first author was supported by an NSF postdoctoral fellowship.
The second author was supported by a Rothschild postdoctoral fellowship.
We thank Sonia Natale who noticed an incompleteness in the proof 
of Lemma 2 in the original version of the paper.

\end{document}